# ACCELERATING DIFFUSIONS[1]

By Chii-Ruey Hwang, Shu-Yin Hwang-Ma and Shuenn-Jyi Sheu

*Academia Sinica, Soochow University and Academia Sinica*

Let $U$ be a given function defined on $\mathbb{R}^d$ and $\pi(x)$ be a density function proportional to $\exp -U(x)$. The following diffusion $X(t)$ is often used to sample from $\pi(x)$,

$$dX(t) = -\nabla U(X(t))\,dt + \sqrt{2}\,dW(t), \qquad X(0) = x_0.$$

To accelerate the convergence, a family of diffusions with $\pi(x)$ as their common equilibrium is considered,

$$dX(t) = (-\nabla U(X(t)) + C(X(t)))\,dt + \sqrt{2}\,dW(t), \qquad X(0) = x_0.$$

Let $L_C$ be the corresponding infinitesimal generator. The spectral gap of $L_C$ in $L^2(\pi)$ ($\lambda(C)$), and the convergence exponent of $X(t)$ to $\pi$ in variational norm ($\rho(C)$), are used to describe the convergence rate, where

$$\lambda(C) = \mathrm{Sup}\{\text{real part of } \mu : \mu \text{ is in the spectrum of } L_C, \mu \text{ is not zero}\},$$

$$\rho(C) = \mathrm{Inf}\Big\{\rho : \int |p(t,x,y) - \pi(y)|\,dy \le g(x)e^{\rho t}\Big\}.$$

Roughly speaking, $L_C$ is a perturbation of the self-adjoint $L_0$ by an antisymmetric operator $C \cdot \nabla$, where $C$ is weighted divergence free. We prove that $\lambda(C) \le \lambda(0)$ and equality holds only in some rare situations. Furthermore, $\rho(C) \le \lambda(C)$ and equality holds for $C = 0$. In other words, adding an extra drift, $C(x)$, accelerates convergence. Related problems are also discussed.

**1. Introduction.** In this paper we prove that by simply adding a weighted divergence-free drift to a reversible diffusion, the convergence to equilibrium is accelerated. In other words, from an algorithmic point of view, the non-reversible algorithm performs better. The analysis is related to the study of antisymmetric perturbations of self-adjoint infinitesimal generators.

Received September 2003; revised May 2004.

[1]Supported in part by NSC Grants NSC 91-2115-M-001-006 and 40222F.

Our investigation is motivated by the following consideration. High- dimensional probability distributions appear frequently in applications. To sample from these distributions directly is not feasible in practice, especially when the corresponding densities are known up to normalizing constants only. One has to resort to approximations. A Markov process with the underlying distribution as its equilibrium is often used to generate an approximation ("MCMC"). How good the approximation is depends on the approximating Markov process and on the specific criterion used for comparison. One may investigate the convergence properties of some particular Monte Carlo Markov processes, or compare the convergence rate within a family of Markov processes (with the same equilibrium) w.r.t. different criteria, or even try to find optimal solutions in that family. Mathematical problems arising from this approach are challenging. Related works may be found in Amit (1991), Amit and Grenander (1991), Frigessi, Hwang and Younes (1992), Frigessi, Hwang, Sheu and di Stefano (1993), Hwang, Hwang-Ma and Sheu (1993), Amit (1996), Athreya, Doss and Sethuraman (1996), Gilks and Roberts (1996), Mengersen and Tweedie (1996), Stramer and Tweedie (1997), Chang and Hwang (1998), Hwang and Sheu (1998, 2000) and Roberts and Rosenthal (2004).

Here we concentrate on the diffusion case. Let $U$ be a given real-valued function defined in $\mathbb{R}^d$ satisfying some smoothness conditions. The underlying distribution $\pi$ is assumed to have a density proportional to $\exp -U(x)$. The following diffusion is commonly used for sampling from its equilibrium $\pi$,

$$(1) \qquad dX(t) = -\nabla U(X(t))\, dt + \sqrt{2}\, dW(t), \qquad X(0) = x_0,$$

where $W(t)$ is the Brownian motion in $\mathbb{R}^d$. For convenience, $\pi$ will be used to denote the underlying probability measure, as well as its probability density. For applications one may consult Grenander and Miller (1994), Miller, Srivastava and Grenander (1995), Srivastava (1996) and references therein.

If a diffusion is regarded as a useful approach to sampling, then it is natural to consider a family of diffusions with $\pi$ as their common equilibrium:

$$(2) \quad dX(t) = -\nabla U(X(t))\, dt + C(X(t))\, dt + \sqrt{2}\, dW(t), \qquad X(0) = x_0,$$

under suitable conditions on $C(x)$. Roughly speaking, the conditions are that $\mathrm{div}(C(x)\exp -U(x)) = 0$ and there is no explosion in (2), that is, $|X(t)|$ does not tend to infinity in a finite time. A strict definition of explosion can be found on page 172 of Ikeda and Watanabe (1989). Exact conditions will be spelled out later. It is easy to pick such a $C$. For example, $C(x) = S(\nabla U(x))$, for any skew symmetric matrix $S$. We are interested in how $C(x)$ influences the convergence of the diffusion (2) to equilibrium.

Hwang, Hwang-Ma and Sheu (1993) focused on a special case, the study of a family of Gaussian diffusions where $2U(x) = (-Dx)\cdot x, -\nabla U(x) =$ $Dx$, $C(x) = SDx$, and where $D$ is a strictly negative-definite real matrix and $S$ is any skew symmetric real matrix. In this case, $\pi(x)$ is Gaussian with mean 0 and covariance matrix $-D^{-1}$ and $X(t)$ is an Ornstein–Uhlenbeck process with drift $(D + SD)x$. Using the rate of convergence of the covariance of $X(t)$ [or together with $EX(t)$] as the criterion, the reversible diffusion with drift $Dx$ (i.e., $C = 0$) is the worst choice and the optimal solution is obtained in this setup.

If $C(x)$ is not zero, then the corresponding diffusion, regarded as a Markov process, is nonreversible. In general, it is difficult to analyze nonreversible processes. We just cite some related works in different settings. In Geman and Geman (1984), Amit and Grenander (1991) and Hwang and Sheu (1998) the convergence properties of some nonreversible Gibb samplers are studied. The ergodicity of systematic sweep in stochastic relaxation, again nonreversible, is investigated in Hwang and Sheu (1992).

Two comparison criteria are considered here. Basic questions such as the acceleration of convergence and the consistency of the comparison w.r.t. these two criteria are answered. Related problems will be discussed in the last section.

Let $\|\cdot\|_p$ and $\|\cdot\|_{p\to q}$ denote the norm in $L^p(\pi)$ and the operator norm from $L^p(\pi)$ to $L^q(\pi)$, respectively, $1 \le p, q \le \infty$. For $p = q = 2$, both norms are simply denoted by $\|\cdot\|$. Let $L_C$ denote the infinitesimal generator of the diffusion $X(t)$ from (2) and, for $C = 0$, let $L = L_0$. Let $T(t) = e^{tL_C}$ denote the corresponding semigroup,

$$T(t)f(x) = E_x f(X(t)) = \int p(t, x, y) f(y)\, dy,$$

where $p(t, x, y)$ is the transition density if it exists. Note that the index $C$ is suppressed from $T(t)$ and $p(t, x, y)$ for the sake of brevity.

We define now the spectral gap of $L_C$ in $L^2(\pi)$ as the first comparison criterion. Since $E_x f(X(t)) \to \pi(f)$ for any starting point $x$, one may consider the average case formulation by averaging the difference $(E_x f(X(t)) - \pi(f))^2$ over the starting point w.r.t. $\pi$:

$$(3) \quad \int (E_x f(X(t)) - \pi(f))^2 \pi(x)\, dx = \|T(t)f - \pi(f)\|^2$$
$$\le \text{constant } \|f - \pi(f)\|^2 e^{2\lambda t},$$

for some $\lambda$ less than or equal to 0, where $\pi(f)$ means integration of $f$ w.r.t. $\pi$. Now consider the worst-case analysis over $f$, then $\|T(t) - \pi\| \le \text{constant } e^{\lambda t}$. The infimum over such $\lambda$'s indicates the convergence rate. This shows that the spectral radius of $T(1)$ in the space $\{f \in L^2(\pi), \pi(f) = 0\}$ is a measure of convergence rate of diffusions to equilibrium. Furthermore, the weak spectral mapping theorem holds between $L_C$ and $e^{tL_C}$ [Nagel (1986), page 91]. Hence, the spectral gap of $L_C$ in $L^2(\pi)$ defined by

(4) $\quad \lambda(C) = \text{Sup}\{\text{real part of } \mu : \mu \text{ in the spectrum of } L_C, \mu \ne 0\}$

is a good candidate to serve as a criterion for the comparison of convergence rates.

The constant in (3) may depend on $C$. If instead we reformulate the inequality in (3) without the constant term,

$$\|T(t)f - \pi(f)\| \leq \|f - \pi(f)\|e^{\lambda t},$$

for some $\lambda$, then the inequality depends only on the behavior of the process around time 0 and the rate will be the same regardless of perturbations [Chen (1992), page 312]. Our interest here is instead in the large-time behavior.

We will always assume that there is no explosion for the diffusions under consideration. Sufficient conditions for nonexplosion may be found, for example, in Proposition 1.10 of Stannat (1999). Since the existence of the transition density is needed in Section 2, for simplicity, we assume that the following assumption holds throughout this paper,

(A1) $\quad C$ and $\nabla U$ are in $L^1(\pi) \cap L^l_{\text{loc}}(\pi) \quad$ for some $l > d$;
$\qquad$ for $f \in C_0^\infty \quad \int (C \cdot \nabla f)\pi = 0.$

Under (A1) there is no explosion in the diffusion (2) and the transition density exists with $\pi$ as its equilibrium distribution [Stannat (1999) and Bogachev, Krylov and Röckner (2001)]. For $f \in C_0^\infty, \int (C \cdot \nabla f)\pi = 0$ means that $C$ is weakly weighted divergence free. This is essential for $\pi$ to be an invariant measure.

Intuitively $L_C$ is a perturbation of a self-adjoint operator $L$ by an antisymmetric operator $C \cdot \nabla$ in $L^2(\pi)$. We are interested in how the spectrum changes. Note that, in general, this perturbation is neither small nor relatively compact. For general references, refer to Kato (1995) and Yosida (1980). $L_C$ is not self-adjoint for nonzero $C$. The spaces considered are real vector spaces of real functions. However, for spectral analysis, one has to consider complex vector spaces. We will make the distinction when it is necessary. Let $C_+$ denote $L_C - L$ and $C_-$ denote $L_{-C} - L$.

We assume that the reversible diffusion (1) w.r.t. $\pi$ has an exponential convergence rate. Equivalently, $L$ has a spectral gap in $L^2(\pi)$, that is,

(A2) $\qquad\qquad\qquad\qquad \lambda(0) < 0.$

The existence of a spectral gap for self-adjoint $L$ has been studied extensively, for example, see Wang (1999).

Under the above two assumptions we prove that $\lambda(C) \leq \lambda(0)$. Furthermore, if $\lambda(0)$ is in the discrete spectrum of $L$, then the equality holds only in some rare situation which is characterized completely. These results are in Theorem 1.

Note that the exponential convergence rate assumption is imposed only on the reversible diffusion. As a consequence of Theorem 1, the perturbed diffusion (2) has a better exponential convergence rate. In other words, adding an extra drift accelerates convergence.

For the nonexplosion of (1), (A2) and $\lambda(0)$ in the discrete spectrum of $L$ to all hold, the following is a sufficient condition [Reed and Simon (1978)]:

$$(5) \qquad 1/2|\nabla U(x)|^2 - \Delta U(x) \to \infty \qquad \text{as } |x| \to \infty.$$

From a probabilistic point of view, one may consider the rate of convergence of $p(t, x, y)$ to $\pi$ in variational norm as a comparison criterion. The variational norm of two probability measures is defined as the supremum of the difference between the two probabilities over all events. This may be regarded as some kind of worst case analysis. Note that the variational norm equals one half of the $L^1(dy)$ distance between the two corresponding densities. Hence, $\rho(C)$ defined below is used as a comparison criterion,

$$(6) \qquad \rho(C) = \text{Inf}\bigg\{\rho \colon \int |p(t, x, y) - \pi(y)|\, dy \le g(x)e^{\rho t}\bigg\}.$$

$g(x)$ may depend on $C$. Usually $g$ is assumed to be essentially locally bounded or locally integrable w.r.t. $\pi$. It needs further study for unrestricted $g$. We prove in Theorems 4 and 5 that $\rho(C) \le \lambda(C)$ and equality holds for the reversible case. Again, using $\rho(C)$ as the comparison criterion, adding an antisymmetric perturbation does help. This result is consistent with the previous one.

It is not clear how the perturbations affect $\rho(C)$ directly. We compare $\rho(C)$ and $\rho(0)$ via $\lambda(C)$ and $\lambda(0)$.

We study the above two criteria only. However, we make the following remarks without giving proofs. Since $T(t)$ is a contractive semigroup in $L^p(\pi)$, for $1 \le p \le \infty$, one may consider (3) in terms of the $L^p$ norm. For a fixed $C$, consider the dependence of the convergence rate on $p$. Note that when (1) is an Ornstein–Uhlenbeck process, $\|T(t) - \pi\|_{1 \to 1}$ does not have exponential convergence rate despite the fact that the corresponding $L$ has a spectral gap in $L^2(\pi)$. For the reversible case, the $L^1(\pi)$ to $L^1(\pi)$ exponential convergence rate is equivalent to the essentially uniform boundedness of $g(x)$ in (6) [Chen (2002)]. If $\|T(t) - \pi\|_{p \to p}$ has exponential convergence rate for some $p \ge 1$ and $\|T(1)\|_{p \to (p+1)}$ is bounded, then $\|T(t) - \pi\|_{q \to q}$ has the same exponential convergence rate for all $q \ge p$.

The use of $\lambda(C)$ as the comparison criterion is studied in Section 2. In Section 3 $\rho(C)$ is the criterion. The relationship between $\lambda(C)$ and $\rho(C)$ is studied. Discussion and related problems are presented in Section 4.

**2. Spectral gap as comparison criterion.** If $\lambda(0)$ is in the discrete spectrum of $L$ in $L^2(\pi)$, then by definition its corresponding eigenspace, denoted by **M**, is finite dimensional. Let $D(\cdot)$ denote "the domain of." Define

$$\varepsilon(f,g) = \int (\nabla f \cdot \nabla g)\pi, \qquad f,g \in C_0^\infty.$$

Then $\varepsilon$ is closable in $L^2(\pi)$. In this section our analysis assumes $\pi(f) = 0, f \in L^2(\pi)$.

THEOREM 1. *If* (A1) *and* (A2) *hold, then* $\lambda(C) \leq \lambda(0)$. *Furthermore, if* $\lambda(0)$ *is in the discrete spectrum of* $L$, *then equality holds if and only if* $C_+$ *or* $C_-$ *leaves a nonzero subspace of* **M** *invariant.*

The following inequality from Stannat [(1999), page 124] will be used repeatedly in the proof.

(7)    If $f \in D(L_C)$    then    $f \in D(\varepsilon)$    and    $\varepsilon(f,f) \leq -\int (L_C f) f \pi$.

We prove first that $\lambda(C) \leq \lambda(0)$. For $f$ with $\|f\| = 1$ and $\pi(f) = 0$, let $g(t) = \|T(t)f\|^2$. $g(0) = 1$ and by (7),

$$g'(t) = 2\int (L_C T(t)f)(T(t)f) \qquad \pi \leq -2\varepsilon(T(t)f, T(t)f) \leq 2\lambda(0)g(t).$$

The above differential inequality implies that the operator norm $\|T(t)\|$ in the space $\{f : f \in L^2(\pi), \pi(f) = 0\}$ is less than or equal to $e^{\lambda(0)t}$. Hence, $\lambda(C) \leq \lambda(0)$.

For a complex valued function $f$, let $f^r$ and $f^i$ denote the real and purely imaginary parts of $f$, respectively.

LEMMA 2. *If* $\lambda(0)$ *is in the discrete spectrum of* $L$, *then there exists a* $\delta > 0$ *such that for any* $a$ *with* $\lambda(0) - \delta \leq a \leq \lambda(0)$ *and any* $b$, $a + ib$ *is not in the continuous spectrum of* $L_C$.

PROOF. Since $\lambda(0)$ is in the discrete spectrum of $L$, there exists $\delta > 0$ such that the spectrum of $L$ restricted to the orthogonal complement of **M** is contained in $(-\infty, \lambda(0) - 2\delta)$.

We prove by contradiction. Assume that there are $a,b$ with $\lambda(0) - \delta \leq a \leq \lambda(0)$ such that $(a+ib)$ is in the continuous spectrum of $L_C$. Let $L_C - (a+ib)$ be denoted by $A$. Then $A$ is one-to-one, the range of $A$ is dense, and $A^{-1}$ is not continuous. To arrive at a contradiction, it suffices to show that for bounded $\{f_n\}$, $Af_n \to 0$ implies $f_n \to 0$.

First we show that $f_n \to 0$ weakly. $Af_n \to 0$, the domain of $A^*$ (the adjoint of $A$) being dense, and the boundedness of $\{f_n\}$ imply the weak convergence of $f_n$ to zero. We claim that $0 \geq \limsup(\varepsilon(f_n, f_n) + a\|f_n\|^2)$. $A(f_n) = ((L_C - a)f_n^r + bf_n^i) + i((L_C - a)f_n^i - bf_n^r))$. Since $A(f_n) \to 0$ and $\{f_n\}$ is bounded, the real part of the inner product of $A(f_n)$ and $f_n$, $\pi(f_n^r(L_C)f_n^r) + \pi(f_n^i(L_C)f_n^i) - a\|f_n\|^2$, goes to zero. By (7), the claim is proved.

Let $f_{n,1}$ be the projection of $f_n$ onto $\mathbf{M}$, $f_{n,2}$ the orthogonal complement. $f_n$ converges weakly, and so do $f_{n,1}$ and $f_{n,2}$. Since $\mathbf{M}$ is finite dimensional, $f_{n,1} \to 0$,

$$0 \geq \limsup(\varepsilon(f_n, f_n) + a\|f_n\|^2) = \limsup(\varepsilon(f_{n,2}, f_{n,2}) + a\|f_{n,2}\|^2)$$
$$\geq \limsup(-\lambda(0) + 2\delta + a)\|f_{n,2}\|^2 \geq \delta \limsup \|f_{n,2}\|^2.$$

Therefore, $f_n \to 0$. $\square$

LEMMA 3. *If $\lambda(C) = \lambda(0)$, then there exists $b$ such that $\lambda(0) + ib$ is in the spectrum of $L_C$.*

PROOF. Let $\{f_n\}$ be a sequence of normalized eigenfunctions of $L_C$ with corresponding eigenvalues $\{a_n + ib_n\}$ such that $a_n < \lambda(0)$ and $a_n \to \lambda(0)$. Then $L_C f_n^r = a_n f_n^r - b_n f_n^i$, $L_C f_n^i = a_n f_n^i + b_n f_n^r$ and, by (7),

$$-a_n = -\pi(f_n^r L_C f_n^r) - \pi(f_n^i L_C f_n^i) \geq \varepsilon(f_n, f_n).$$

As in the last part of the proof of Lemma 2,

$$\lambda(0) - a_n \geq \varepsilon(f_n, f_n) + \lambda(0) = \varepsilon(f_{n,2}, f_{n,2}) + \lambda(0)\|f_{n,2}\|^2 \geq \delta\|f_{n,2}\|^2,$$

where $f_{n,2}$ and $\delta$ are as in Lemma 2. Hence, $f_{n,2} \to 0$. Since the projection $\{f_{n,1}\}$ of $\{f_n\}$ onto the finite-dimensional $\mathbf{M}$ is bounded, there exists a convergent subsequence of $\{f_{n,1}\}$. For convenience, the same index $n$ will be used. We have $f_n$ converging to some $f$ in $\mathbf{M}$. Note that the spectral mapping theorem holds for point spectrums. Hence,

$$e^{a_n + ib_n} f_n = T(1)f_n \to T(1)f \quad \text{and} \quad e^{ib_n} \to e^{-\lambda(0)}\pi(fT(1)f).$$

Therefore, there exists some $b$ such that $\lambda(0) + ib$ and $e^{\lambda(0)+ib}$ are eigenvalues of $L_C$ and $T(1)$ with the same eigenfunction $f$. If $\lambda(0)$ is a limit point of the real parts of the residual spectrum of $L_C$, then we can repeat the above proof for the adjoint of $L_C$ which is $L_{-C}$. Hence, there exists some $b$ such that $\lambda(0) + ib$ is in the point spectrum or the residual spectrum of $L_C$. Since there is no continuous spectrum in the neighborhood of $\lambda(0)$, this completes the proof. $\square$

PROOF OF THEOREM 1. If $\lambda(0)$ is the real part of an eigenvalue of $L_C$ with a normalized eigenfunction $f + ig$, then by (7) and the definition of the Dirichlet form $\varepsilon$,

$$-\lambda(0) \geq \varepsilon(f, f) + \varepsilon(g, g) \geq -\lambda(0).$$

Then $f$ and $g$ are in **M** and $C_+$ maps the subspace spanned by $f$ and $g$ into itself. If for some $b$, $\lambda(0) + ib$ is in the residual spectrum of $L_C$, then $\lambda(0) - ib$ is an eigenvalue of the adjoint operator $L_{-C}$. Hence, $C_-$ leaves a nonzero subspace of **M** invariant.

The proof of the other direction is obvious. □

REMARK. It seems that a stronger result should hold: if $\lambda(C) = \lambda(0)$, then $\lambda(0)$ is the real part of an eigenvalue of $L_C$. If this is the case, Theorem 1 has a stronger form: the equality holds iff $C_+$ leaves a nonzero subspace of **M** invariant. If (5) holds, then $(L-a)^{-1}$ is compact for $a$ in the resolvent of $L$ [Reed and Simon (1978)]. And the stronger statements hold.

REMARK. As mentioned in the Introduction, the existence of the transition density is not needed here. A weaker assumption than (A1) suffices, for example, $C$ and $\nabla U$ are in $L^1(\pi) \cap L^2_{\text{loc}}(\pi)$ [Stannat (1999)].

**3. Convergence rate in variational norm as criterion.** Under (A1) the transition density $p(t, x, y)$ exists. Let $p_t(x, y)$ denote $p(t, x, y)/\pi(y)$; $p_t(x, y)$ is locally Hölder [Bogachev, Krylov and Röckner (2001)].

THEOREM 4. *In addition to* (A1) *and* (A2), *if* $\nabla U$ *and* $C$ *are locally bounded, then there exists a locally bounded function* $g$ *such that*
$$\int |p_t(x,y) - 1|\pi(y)\,dy \leq g(x)e^{\rho(c)t}.$$
*Moreover,* $\rho(C) \leq \lambda(C)$.

PROOF.
$$\int |p_t(x,y) - 1|\pi(y)\,dy$$
$$= \int \left|\int (p_1(x,z)p_{t-1}(z,y) - 1)\pi(z)\,dz\right|\pi(y)\,dy$$
$$= \int \left|\int (p_1(x,z)p^*_{t-1}(y,z) - 1)\pi(z)\,dz\right|\pi(y)\,dy$$

(* denotes the adjoint process)

$$= \int \left|\int p^*_{t-1}(y,z)p_1(x,z)\pi(z)\,dz - \int p_1(x,z)\pi(z)\,dz\right|\pi(y)\,dy$$
$$= \int |T^*(t-1)(p_1(x,\cdot))(y) - \pi(p_1(x,\cdot))|\pi(y)\,dy$$
$$\leq \left(\int |T^*(t-1)(p_1(x,\cdot))(y) - \pi(p_1(x,\cdot))|^2 \pi(y)\,dy\right)^{1/2}$$
$$\leq \text{constant}\|p_1(x,\cdot) - 1\|e^{\lambda(C)t}.$$

The last inequality holds if $p_1(x,\cdot)$ is in $L^2(\pi)$.

We now claim that $\int p_1^2(x,z)\pi(z)\,dz$ is locally bounded. Since $\nabla U$ and $C$ are locally bounded, by a local Harnack inequality [Theorem 1.1 in Trudinger (1968)],

$$\begin{aligned}&\forall x \text{ in } \mathbb{R}^d, \forall N > 0, \forall f \text{ with } \pi(f) = 1 \text{ and } f \geq 0,\\ &\quad\sup_{y\in B(x,N/2)} T(s)f(y) \leq C(N,x) \inf_{y\in B(x,N/2)} T(2s)f(y),\end{aligned} \quad (8)$$

where the constant $C(N,x)$ depends only on $N$ and $x$, and $B(x,N/2)$ denotes a ball in $\mathbb{R}^d$ with center $x$ and radius $N/2$. For $y$ and $z$ in $B(x,N/2)$,

$$\begin{aligned}\int p_s(z,u)f(u)\pi(u)\,du &= T(s)f(z) = \frac{\int_{B(x,N/2)} T(s)f(z)\pi(y)\,dy}{\pi(B(x,N/2))}\\ &\leq \frac{C(N,x)\int_{B(x,N/2)} T(2s)f(y)\pi(y)\,dy}{\pi(B(x,N/2))} \leq \frac{C(N,x)}{\pi(B(x,N/2))};\end{aligned}$$

for $f$ satisfying (8), we have

$$\sup_y p_s(z,y) \leq \frac{C(N,x)}{\pi(B(x,N/2))}.$$

Now let $g(x) = \sup_y p_s(x,y)$, then $g$ is locally bounded and

$$\int p_s^2(x,y)\pi(y)\,dy \leq g^2(x).$$

This also establishes that $\rho(C) \leq \lambda(C)$. $\square$

REMARK. The local boundedness assumption in Theorem 4 is not needed for the reversible case, since $\int p_1^2(x,y)\pi(y)\,dy = p_2(x,x)$ is locally bounded.

The following theorem implies that for the reversible case, $\rho(0) = \lambda(0)$.

THEOREM 5. *For the reversible case, if there exists some $g$ in $L^1_{\text{loc}}(\pi)$ such that*

$$\int |p_t(x,y) - 1|\pi(y)\,dy \leq g(x)e^{\rho t} \qquad \text{then } \|T(t) - \pi\| \leq e^{\rho t}.$$

PROOF. For the reversible case, $T(t)$ is self-adjoint in $L^2(\pi)$,

$$\|T(t)f\|^2 = \pi(fT(2t)f).$$

For $f$ with $\pi(f) = 0$, $f \in C^\infty$, $f = c_0$ outside $B_N$, where $B_N$ denotes a ball in $\mathbb{R}^d$ centered at 0 with radius $N$ and $c_0$ a constant,

$$\|T(t)f\|^2 = \pi(fT(2t)f) = \pi((f-c_0)T(2t)(f-c_0)) - c_0^2$$

$$= \int_{B_N} (f(x) - c_0) \left( \int_{B_N} (p_{2t}(x,y) - 1)(f(y) - c_0)\pi(y)\,dy \right) \pi(x)\,dx$$

$$\leq \|f - c_0\|_\infty^2 \int_{B_N} \int_{B_N} |p_{2t}(x,y) - 1| \pi(y)\,dy\,\pi(x)\,dx$$

$$\leq \|f - c_0\|_\infty^2 \int_{B_N} g(x) e^{2\rho t} \pi(x)\,dx \leq C(N,f) e^{2\rho t}.$$

By Lemma 2.2 in Röckner and Wang (2001), for $s \leq t$ and $\pi(f^2) = 1$,

$$\|T(s)f\|^2 \leq (\|T(t)f\|^2)^{s/t} \leq C(N,f)^{s/t} e^{2\rho s}.$$

The equalities hold at $s = 0$. Now take a derivative w.r.t. $s$ and evaluate at 0. Then

$$-2\varepsilon(f,f) \leq 1/t \log C(N,f) + 2\rho.$$

Letting $t \to \infty$, we have $\varepsilon(f,f) \geq -\rho$.

For any $h \in C_0^\infty$,

$$\text{let } f = \frac{h - \pi(h)}{\|h - \pi(h)\|} \qquad \text{then } \varepsilon(f,f) \geq -\rho,$$

and $\pi(h^2) \leq \frac{1}{-\rho}\varepsilon(h,h) + \pi^2(h)$. Hence, we have proven $\lambda(0) \leq \rho$. $\square$

**4. Discussion and related problems.** Our theorems give only general and qualitative information. The proofs do not reveal how the rate of convergence depends on $C$. Intuitively, multiplying $C$ by a large $k$ should speed up convergence. However, examples in Hwang, Hwang-Ma and Sheu (1993) show the contrary. It is not clear which part of $C$ contributes to acceleration. Most of the questions discussed below are based on $\lambda(C)$. Similar questions can be formulated for $\rho(C)$.

Now consider families of diffusions (algorithms) defined by (2) with index $C$ satisfying various conditions. What is the best algorithm within a certain family? For example, let **G** denote the family of diffusions with $C$ satisfying the general conditions described in the previous sections and **S** the family of diffusions with $C = S(\nabla U)$ for any skew symmetric matrix $S$, respectively. One may ask for the optimal values and minimizers in the following two problems:

($\lambda$1) $\text{Inf}_{C \in \mathbf{G}} \lambda(C)$.
($\lambda$2) $\text{Inf}_{C \in \mathbf{S}} \lambda(C)$.

For ($\lambda$1), even the simple question "Is $\text{Inf}_{C \in \mathbf{G}} \lambda(C) = -\infty$?" remains unanswered. For Gaussian diffusions, the optimal structure of ($\lambda$2) is known [Hwang, Hwang-Ma and Sheu (1993)]. Note that for the Gaussian case, the perturbation $C$ in ($\lambda$2) is linear. However, Hwang and Sheu (2000) showed that a quadratic $C$ has a better rate. Problem ($\lambda 2$) remains open for the general case.

Basically, the problem is to find the best "spectral gap" in a family of elliptic operators. Similar questions may be discussed on compact Riemannian manifolds. We consider the following generic case on the two-dimensional torus: let $L_C = \Delta + C \cdot \nabla$ with divergence free $C$,

($\lambda 3$) $\inf_C \lambda(C)$.

Again, what is the best solution and is it finite?

Obviously, $\lambda(0) < 0$ implies $\lambda(C) < 0$, but how about the other way around? That is, if $L_C$ has a spectral gap, does $L$? If the answer is negative, then perturbations can drastically change fundamental convergence properties. We proved that $\rho(C) \leq \lambda(C)$, but when does equality hold? If there is no spectral gap for $L$, how does the antisymmetric perturbation accelerate convergence?

**Acknowledgments.** C.-R. Hwang would like to thank Professor M. F. Chen and Professor F. Y. Wang for many inspiring discussions. Part of the work was written while C.-R. Hwang was visiting the CMLA, ENS-Cachan, France; he would like to thank Professor L. Younes for suggestions and hospitality.

C.-R. Hwang  
S.-J. Sheu  
Institute of Mathematics  
Academia Sinica  
Taipei  
Taiwan 11529  
e-mail: crhwang@sinica.edu.tw

S.-Y. Hwang-Ma  
Department of Business Mathematics  
Soochow University  
Taipei  
Taiwan 10001